\documentclass[12pt]{amsart}

\newtheorem{theorem}{Theorem}[section]
\newtheorem{lemma}[theorem]{Lemma}
\newtheorem{proposition}[theorem]{Proposition}

\theoremstyle{definition}   

\newtheorem{example}[theorem]{Example}

\theoremstyle{remark}

\numberwithin{equation}{section}

\usepackage{times}
\usepackage{enumerate}
\usepackage{tfrupee}
\usepackage{tikz}
\usetikzlibrary{chains,fit}
\usetikzlibrary{shapes,snakes}
\usetikzlibrary{graphs}
\usepackage{graphicx,adjustbox}
\usepackage{hyperref}
\usepackage{amsmath,amssymb}

\usepackage{amscd}
\usepackage{graphicx}
\usepackage[all]{xy}

\title[Moh's example of algebroid space curves]
{Moh's example of algebroid space curves}
\author{
Ranjana Mehta
\and
Joydip Saha
\and
Indranath Sengupta
}
\date{}

\address{\small \rm  Discipline of Mathematics, IIT Gandhinagar, Palaj, Gandhinagar, 
Gujarat 382355, INDIA.}
\email{ranjssj16@gmail.com}

\address{\small \rm  Discipline of Mathematics, IIT Gandhinagar, Palaj, Gandhinagar, 
Gujarat 382355, INDIA.} 
\email{saha.joydip56@gmail.com}
\thanks{The second author thanks SERB, Government of India for the post-doctoral fellow 
position at IIT Gandhinagar, through the research project EMR/2015/000776.}

\address{\small \rm  Discipline of Mathematics, IIT Gandhinagar, Palaj, Gandhinagar, 
Gujarat 382355, INDIA.}
\email{indranathsg@iitgn.ac.in}
\thanks{The second author is the corresponding author; supported by the 
MATRICS research grant MTR/2018/000420, sponsored by the SERB, Government of India.}

\date{}



\begin{document}
\begin{abstract}
In this paper we revisit the family of algebroid space curves 
defined by Moh and find an explicit minimal generating set for 
the defining ideal.
\end{abstract}

\maketitle

Moh \cite{moh1} defined the following family of algebroid space curves to 
exhibit that its defining ideal requires an arbitrarily large number of 
generators. Let $k$ be a field of characteristic zero, $n$ an odd positive 
integer, $m=(n+1)/2$ and $\lambda$ an integer coprime with $m$ such that 
$\lambda > n(n+1)m$. Let $\rho_{n}: k[[x,y,z]] \rightarrow k[[t]]$ be a 
mapping defined by 
$$\rho_{n}(x)= t^{nm}+t^{nm+\lambda},\, \rho_{n}(y)= t^{(n+1)m},\,\rho_{n}(z)= t^{(n+2)m}.$$
We define $\mathfrak{p}_{n} = \ker(\rho_{n})$ and $P_{n}=\mathfrak{p}_{n}\cap k[x,y,z]$. 
Moh proved in \cite{moh2} that the minimal number of generators of $\mathfrak{p}_{n}$ 
is $\mu(\mathfrak{p}_{n})= n+1$. Moh also proved that the ideal 
$\mathfrak{p}_{n}$ has a determinantal structure. However, it is not easy to 
read the generators explicitly from Moh's work. This paper is devoted to 
understanding Moh's results and use them to write an explicit minimal 
generating set. We have given a very precise algorithm to compute the 
polynomials generating the ideal $\mathfrak{p}_{n}$ minimally. However, 
it is not very clear to us whether they are the same as those conceived 
in Moh's work or whether there is any underlying determinantal structure. 
\medskip

The generators in explicit 
form often become important especially when one is trying to understand 
issues like smoothness of blowups for such curves, which is particularly 
interesting in the case $n=1$, when it is known to be a complete intersection. 
We have tried out two different approaches for constructing a minimal 
generating set. For the special case of $n=1$, our 
approach would be to find a Gr\"{o}bner basis first for the contracted 
ideal $P_{1}=\mathfrak{p}_{1}\cap k[x,y,z]$ in the polynomial ring 
$k[x,y,z]$ in section 2 and use that to construct a minimal set of generators for 
$\mathfrak{p}_{1}$ in the power series ring in section 3. However, 
this method does not prove to be effective for higher values of $n$ 
because of two reasons. For a general $n$, it is not easy to compute a 
Gr\"{o}bner basis for $P_{n}$ in the polynomial ring $k[x,y,z]$ and also 
the difficulty arises in extracting a minimal generating set for 
$\mathfrak{p}_{n}$ in the power series ring $k[[x,y,z]]$. The 
case for general $n$ has been worked out in detailed in section 4 and we have 
adopted a more direct method where we have used Moh's results. Starting 
with two generators for the ideal $\mathfrak{p}_{n}$, with the help of initial 
guesses through computations with \cite{singular}, we have built the other 
generators iteratively. The method is described stepwise in section 4 and also 
illustrated through an example.

\section{Moh's techniques}
Let us first recall the technique used in \cite{moh1} and \cite{moh2}, 
where he defined the map $\sigma: k[[x,y,z]] \rightarrow k[[x,y,z]]$ as $\sigma(x) = x^{n}$, 
$\sigma(y)=y^{n+1}$, $\sigma(z) = z^{n+2}$. The $\sigma$-order of $f(x,y,z)\in k[[x,y,z]]$ 
is $\mbox{ord}(\sigma (f(x,y,z)))$, the order of the power series $\sigma (f(x,y,z))$. 
The $\sigma$ leading form of $f(x,y,z)$ is $\sigma ^{-1}$ (the leading form of 
$\sigma (f(x,y,z))$) and $f(x,y,z)\in k[[x,y,z]]$ is said to be $\sigma$-homogeneous 
if $f(x,y,z)$= the $\sigma$-leading form of $f(x,y,z)$. Therefore, 
$f$ is $\sigma$-homogeneous if and only if $\sigma(f)$ is homogeneous. 
We use the notation $f^{\sigma}$ to write the $\sigma$-leading 
form of $f$ and the notation $f^{\tau} = f - f^{\sigma}$ to denote the 
\textit{tail} of $f$. In section 4, where we discuss the situation for 
general $n$, which would heavily depend on description of $f^{\sigma}$ 
and $f^{\tau}$. In fact, our technique would be to create $n+1$ polynomials 
in $\mathfrak{p}_{n}$ by creating the $\sigma$-leading forms and the tails 
separately. Finally, Moh's theorem will be used to prove that these 
polynomials indeed generate $\mathfrak{p}_{n}$ minimally. Let us recall the 
main theorem proved by Moh in \cite{moh2}, which is frequently used in our argument. 
The theorem requires the definition of the vector space $V_{r}$ consisting 
of the $\sigma$-leading forms of elements in $\mathfrak{p}_{n}$ of $\sigma$-order 
$r$, since the $\sigma$-leading forms are $\sigma$-homogeneous. 
\medskip
 
\noindent\textbf{Theorem (Moh; \cite{moh2})} \label{gen}
The prime ideal $\mathfrak{p}_{n}$ needs at least $n+1$ generators. There are 
$f_{1},\ldots,f_{n+1}\in \mathfrak{p}_{n}$ such that:
\begin{enumerate}
\item the $\sigma$-leading form of $f_{i}$ belong to the set $V_{n^2+n+i-1}$, for $i=1,\ldots, n$, where 
\begin{eqnarray*}
V_{r} & = &\{\sigma{\rm -homogeneous \, form \, of}\, \sigma{\rm -order}\, r\} \cap \\
{} & {} & \{\sigma{\rm -leading \, forms \, of \, elements \, in}\, \mathfrak{p}_{n}\}\cup \textbf{0};
\end{eqnarray*}
\item $x$ (the $\sigma$-leading form of $f_{1}$) and the $\sigma$-leading form of $f_{n+1}$ generate 
$V_{n^2+2n}$.
\end{enumerate}
Moreover, any $f_{1},\ldots,f_{n+1}\in \mathfrak{p}_{n}$ satisfying conditions (1) and (2) generate $\mathfrak{p}_{n}$.

\section{Gr\"{o}bner basis of the ideal $P_{1}$ in $k[x,y,z]$}
Our motivation behind this study is to understand Moh's class of curves from the 
computational perspective. We will see that it would be helpful if we rewrite Moh's 
examples by introducing some new parameters. Recall that $n$ is an odd positive 
integer, $m=(n+1)/2$ and $\lambda$ is an integer such that $\lambda > n(n+1)m$ with 
$\gcd(\lambda, m)=1$. Let us write $\lambda = A+Bp+r$, where 
$A= n(n+1)m$; $B= (n+2)m$; $1\leq r\leq B \, \text{and}\, p\in \mathbb{Z}_{\geq0}$. 
Note that $\gcd(\lambda, m)=1$ implies  $\gcd(r, m)=1$. In the case of $n=1$, 
we have $m=1$, $\lambda>2$, $A=2$, $B=3$, $r\in\{1,\, 2,\, 3\}$ and 
$p\in \mathbb{Z}_{\geq 0}$. 
\medskip

\begin{proposition}\label{gbth} 
Let us consider the lexicographic monomial order induced by $x>y>z$ on $k[x,y,z]$. Let us write 
\begin{itemize}
\item $g_{1} = y^{3}- z^{2}$,
\item $g_{2} = xz-y^{2}-y^{(3-r)}z^{(r+p)}$,
\item $g_{3} = xy-z-y^{(4-r)}z^{(r+p-1)}$,
\item $g_{4} = x^{2}-y-xy^{(3-r)}z^{(r+p-1)}-y^{\frac{1}{2}(r-2)(3r-5)}z^{(-r^2+4r+p-2)}$;
\end{itemize}
where $r$ and $p$ are defined above. The set $\mathcal{G}=\{g_{1},g_{2},g_{3},g_{4}\}$ 
forms a Gr\"{o}bner basis for the ideal $I=\langle g_{1},g_{2},g_{3},g_{4}\rangle$, with 
respect to the chosen monomial order.
\end{proposition}

\proof  We proceed by Buchberger's criterion \cite{coxetal} and examine each S-polynomial separately.
\begin{enumerate}[(i)]
\item We have $\gcd(\mathrm{LT}(g_{1}),\mathrm{LT}(g_{2}))=1$ and $\gcd(\mathrm{LT}(g_{1}),\mathrm{LT}(g_{4}))=1$. 
Therefore $S(g_{1},g_{2})\longrightarrow_{\mathcal{G}} 0$ and $S(g_{1},g_{4})\longrightarrow_{\mathcal{G}} 0$.
\item We have $S(g_{1},g_{3})= -xz^{2}+y^{2}z+y^{(6-r)}z^{(r+p-1)}$. Therefore, 
$S(g_{1},g_{3})=-z\cdot g_{2} +y^{(3-r)}z^{(r+p-1)}\cdot g_{1}\longrightarrow_{\mathcal{G}} 0$.
\item We have $S(g_{2},g_{3})= -g_{1}\longrightarrow_{\mathcal{G}} 0 $.
\item We have $S(g_{2},g_{4})= -xy^{2}+yz+y^{\frac{1}{2}(r-2)(3r-5) }z^{(-r^2+4r+p-1)}$.  If $r$ is either $1$ or $2$, 
then $S(g_{2},g_{4})=-y^{(2-r)}z^{(r-1+p)}\cdot g_{1}- y\cdot g_{3}\longrightarrow_{\mathcal{G}} 0$. If $r=3$, then 
$S(g_{2},g_{4})= -y\cdot g_{3}\longrightarrow_{\mathcal{G}} 0$.
\item We have $S(g_{3},g_{4})=-xz+y^{2}+y^{(\frac{1}{2}r^2-\frac{11}{2}r+6)}z^{(-r^2+4r-2+p)}$. 
If $r$ is either $1$ or $2$, then $ S(g_{3},g_{4})=-(g_{2})\longrightarrow_{\mathcal{G}} 0$. If $r=3$, then 
$S(g_{3},g_{4})=-g_{2}+z^{(1+p)}\cdot {g_1}\longrightarrow_{\mathcal{G}} 0$. \qed
\end{enumerate}
\medskip

\begin{theorem}\label{prime}
The ideal $I=\langle g_{1},g_{2},g_{3},g_{4}\rangle$ is a prime ideal in the polynomial ring $k[x,y,z]$.
\end{theorem}

\proof Let us consider the lexicographic monomial order in $k[x,y,z]$ induced by $x>y>z$. 
By Proposition \ref{gbth} the set $\mathcal{G}=\{g_{1},g_{2},g_{3},g_{4}\}$ forms a 
Gr\"{o}bner basis 
for the ideal $I$ with respect to the said order. Suppose $fh\in I$ and $f\notin I$, 
$h\notin I$. Without loss of generality, we may assume that no term of $f$ and $h$ 
is divisible by the leading terms $\mathrm{LT}(g_{i})$, $1\leq i\leq 4$, as we may 
replace them by their remainders. Therefore, in particular, neither $\mathrm{LT}(f)$ 
nor $\mathrm{LT}(h)$ is divisible by any one of those leading terms. We know that 
$\mathrm{LT}(g_{1})=y^{3}$, 
$\mathrm{LT}(g_{2})=xy$, $\mathrm{LT}(g_{3})=xz$ and $\mathrm{LT}(g_{4})=x^{2}$. We 
claim that, $\mathrm{LT}(g_{i})$, $2\leq i\leq 4$ can not divide any monomial in 
$\mathrm{supp}(f)\cup\mathrm{supp}(h)$, otherwise, we use division algorithm to 
reduce $f$ (or $h$) modulo $\mathcal{G}=\{g_{1},g_{2},g_{3},g_{4}\}$ to the nonzero 
remainder $R$ and work with $R$ instead of $f$ (or $h$). 
\medskip

\noindent\textbf{Case 1.}\, Let $y^{3}\mid \mathrm{LT}(fh)=\mathrm{LT}(f)\mathrm{LT}(h)$. 
We assume that $x\nmid \mathrm{LT}(fh)$ otherwise we consider case 2.  Note that the 
indeterminate $x$ is the largest among $x, y, z$ and it does not divide $\mathrm{LT}(f)$; 
therefore $x$ does not divide any monomial in the support of $f$. 
Since $y^{3}\nmid \mathrm{LT}(f) $ and $y^{3}\nmid \mathrm{LT}(h)$, therefore, 
$y^{3}$ does not divide the other terms of $f$ and $h$. We may write 
$$f=y^{2}p_{1}(z)+yp_{2}(z)+p_{3}(z)\in k[y,z];$$ where $p_{i}(z)\in k[z]$ for 
$1\leq i\leq 3$. Similarly we can write 
$$h=y^{2}q_{1}(z)+yq_{2}(z)+q_{3}(z)\in k[y,z];$$ 
where $q_{i}(z)\in k[z]$ for $1\leq i\leq 3$. Therefore 
\begin{eqnarray*}
fh & = & y^{4}p_{1}(z)q_{1}(z)+y^{3}(p_{1}(z)q_{2}(z)+p_{2}(z)q_{1}(z))+ y^{2}(p_{1}(z)q_{3}(z)\\
{} & {} & +p_{2}(z)q_{2}(z)+p_{3}(z)q_{1}(z))+y(p_{2}(z)q_{3}(z)+p_{3}(z)q_{2}(z))+p_{3}(z)q_{3}(z).
\end{eqnarray*}
\noindent We divide by $g_{1}$ to get,
\begin{align*}
& yz^{2}p_{1}(z)q_{1}(z)+z^{2}\left[p_{1}(z)q_{2}(z)+p_{2}(z)q_{1}(z)\right]
   + y^{2}[p_{1}(z)q_{3}(z)+p_{2}(z)q_{2}(z)\\
&+p_{3}(z)q_{1}(z)]+y[p_{2}(z)q_{3}(z)+p_{3}(z)q_{2}(z)]+p_{3}(z)q_{3}(z)\in I.
\end{align*}
This leads to a contradiction because the leading term of the above expression 
is not divisible by any $\mathrm{LT}(g_{i})$, $1\leq i\leq 4$.
\medskip
   
\noindent\textbf{Case 2.} Let $xy\mid\mathrm{LT}(fh)=\mathrm{LT}(f)\mathrm{LT}(h)$. 
Without loss of generality we may assume that $x\mid \mathrm{LT}(f)$ and 
$y\mid \mathrm{LT}(h)$. Since $\mathrm{LT}(f)$ and $\mathrm{LT}(h)$ are 
not divisible by any $\mathrm{LT}(g_{i})$, $1\leq i\leq 4$, we have $x^{2}, y,z$ 
not dividing $\mathrm{LT}(f)$ and $x,y^{3}$ not dividing $\mathrm{LT}(h)$. 
Hence $x$ does not divide any monomial in $\mathrm{supp}(h)$. Furthermore, 
$y^{3}$ does not divide any monomial in $\mathrm{supp}(f)\cup \mathrm{supp}(h)$. 
Therefore we have $f=x+p(y,z)$, where 
$$p(y,z)=y^{2}p_{1}(z)+yp_{2}(z)+p_{3}(z)\in k[y,z]$$ 
and 
$$h=y^{2}q_{1}(z)+yq_{2}(z)+q_{3}(z),$$ 
where $p_{i}(z), q_{i}(z)\in k[z]$ for $1\leq i\leq 3$. Therefore 
\begin{eqnarray*}
fh & = & xy^{2}q_{1}(z)+xyq_{2}(z)+xq_{3}(z)+y^{2}q_{1}(z)p(y,z)\\
{} & {} & +yq_{2}(z)p(y,z)+q_{3}(z)p(y,z)\in I;
\end{eqnarray*} 
and dividing by $g_{3}$ we get
\begin{align}\label{ex1}
&[z+y^{(4-r)}z^{(r+p-1)}][yq_{1}(z)+ q_{2}(z)]+xq_{3}(z)+y^{2}q_{1}(z)p(y,z)\\ 
&+yq_{2}(z)p(y,z)+q_{3}(z)p(y,z)\in I . \nonumber
\end{align}
 
\noindent Let $q_{3}(z)=zr(z)+c$, then, from expression \ref{ex1} we have
\begin{align*}
&[z+y^{(4-r)}z^{(r+p-1)}][yq_{1}(z)+ q_{2}(z)]+x[zr(z)+c]+y^{2}q_{1}(z)p(y,z)\\
&+yq_{2}(z)p(y,z)+[zr(z)+c]p(y,z)\in I .
\end{align*}
Dividing the above expression by $g_{2}$ we get 
\begin{align*}
&[z+y^{(4-r)}z^{(r+p-1)}][yq_{1}(z)+ q_{2}(z)]+[y^{2}+y^{(3-r)}z^{(r+p)}]r(z)+cx\\
&+y^{2}q_{1}(z)p(y,z)+yq_{2}(z)p(y,z)+[zr(z)+c]p(y,z)\in I .
\end{align*}
  
\noindent If $c\neq 0$ then the leading term of the above expression is $cx$, 
which is not divisible by any $\mathrm{LT}(g_{i})$, $1\leq i\leq 4$, and this leads 
to a contradiction. If $c=0$, then substituting $p(y,z)$ by 
$y^{2}p_{1}(z)+yp_{2}(z)+p_{3}(z)$ and dividing by $g_{1}$ we get 
\begin{align*}
& [z+y^{(4-r)}z^{(r+p-1)}][yq_{1}(z)+ q_{2}(z)]+[y^{2}+y^{(3-r)}z^{(r+p)}]r(z)\\
& +y^{2}[q_{2}(z)p_{2}(z)+zr(z)p_{1}(z)]\\
& +y[z^{2}q_{1}(z)p_{1}(z)+q_{2}(z)p_{3}(z)+zr(z)+p_{2}(z)]\\
& +z^{2}[q_{1}(z)p_{2}(z)+q_{2}(z)p_{1}(z)]+zr(z)p_{3}(z)\in I
\end{align*}
Again we consider cases depending on different values of $r$. If $r=3$, 
then the leading term of the above expression is not 
divisible by any $\mathrm{LT}(g_{i})$, $1\leq i\leq 4$, which is impossible.
If $r=1$, then division by $g_{1}$ gives us
\begin{align*}
& [z+z^{(2+p)}][yq_{1}(z)+ q_{2}(z)]+[y^{2}+y^{2}z^{(1+p)}]r(z)\\
&+y^{2}[q_{2}(z)p_{2}(z)+zr(z)p_{1}(z)]\\
& +y[z^{2}q_{1}(z)p_{1}(z)+q_{2}(z)p_{3}(z)+zr(z)+p_{2}(z)]\\
& +z^{2}[q_{1}(z)p_{2}(z)+q_{2}(z)p_{1}(z)]+zr(z)p_{3}(z)\in I
\end{align*}
Again leading term of the above expression is not divisible by any $\mathrm{LT}(g_{i})$, 
$1\leq i\leq 4$, which gives a contradiction. Finally, if $r=2$, we divide the expression 
by $g_{1}$ and get
\begin{align*}
& z[yq_{1}(z)+ q_{2}(z)]+z^{(p+3)}q_{1}(z)+y^{2}z^{(p+1)} q_{2}(z)+[y^{2}+yz^{(2+p)}]r(z)\\
& +y^{2}[q_{2}(z)p_{2}(z)+zr(z)p_{1}(z)]\\
& +y[z^{2}q_{1}(z)p_{1}(z)+q_{2}(z)p_{3}(z)+zr(z)+p_{2}(z)]\\
& +z^{2}[q_{1}(z)p_{2}(z)+q_{2}(z)p_{1}(z)]+zr(z)p_{3}(z)\in I
\end{align*}
Leading term of the above expression is not divisible by any $\mathrm{LT}(g_{i})$, 
$1\leq i\leq 4$, which gives a contradiction.
\medskip
 
\noindent\textbf{Case 3.} Let $xz\mid\mathrm{LT}(fh)=\mathrm{LT}(f)\mathrm{LT}(h)$. 
We assume that $y\nmid \mathrm{LT}(fh)$, otherwise $xy\mid \mathrm{LT}(fh)$ and we are
back to the Case 2. Without loss of generality we assume that $x\mid \mathrm{LT}(f)$ 
and $z\mid \mathrm{LT}(h)$. Since $x^{2}, y,z$ do not divide $\mathrm{LT}(f)$ and 
$x,y$ do not divide $\mathrm{LT}(h)$, we have $f=x+y^{2}p_{1}(z)+yp_{2}(z)+p_{3}(z)$ 
and $h=q(z)$, where $q(z),p_{i}(z)\in k[z]$, $1\leq i\leq 3$. Then, 
$$fh=xq(z)+[y^{2}p_{1}(z)+yp_{2}(z)+p_{3}(z)]q(z).$$ Suppose that 
$q(z)=zt(z)+d$; after dividing by $g_{2}$ we get
$$ [y^{2}+y^{(3-r)}z^{(r+p)}]t(z)+dx+[y^{2}p_{1}(z)+yp_{2}(z)+p_{3}(z)]q(z)\in I.$$ If $d\neq 0$ 
then the leading term of the above expression is $dx$, which is not divisible by any 
$\mathrm{LT}(g_{i})$, $1\leq i\leq 4$ -- a contradiction. If $d=0$, then also the leading 
term of the above expression is not divisible by any $\mathrm{LT}(g_{i})$, $1\leq i\leq 4$ 
-- a contradiction.
\medskip

\noindent\textbf{Case 4.} Let $x^{2}\mid \mathrm{LT}(fh)$; then $x$ divides $\mathrm{LT}(f)$ 
and $\mathrm{LT}(h)$ both. Again by the same argument as in the previous cases we have 
$$f=x+y^{2}p_{1}(z)+yp_{2}(z)+p_{3}(z)$$ 
and 
$$h=x+y^{2}q_{1}(z)+yq_{2}(z)+q_{3}(z),$$ 
where $p_{i},q_{i}(z)\in k[z]$, for $1\leq i\leq 3$. We argue similarly and 
keep on dividing by $g_{i}$'s finally arriving at a contradiction.\qed
\medskip

\begin{proposition}\label{irreducible}
The polynomials $g_{3}$ and $g_{4}$ are irreducible in the ring $k[x,y,z]$ 
and $k[[x,y,z]]$. Hence $\{g_{4},g_{3}\}$ as well as $\{g_{3},g_{4}\}$ form 
regular sequences of length $2$.
\end{proposition}

\proof We write $g_{3}=yx+b(y,z)$, where $b(y,z)=-[z+y^{(4-r)}z^{(r+p-1)}]$. 
Note that $z\mid b(y,z)$, $z\nmid y$ and $z^{2}\nmid b(y,z)$. We can use 
Eisenstein criterion with respect to the prime element $z$ and 
prove that $g_{3}$ is irreducible in $k[x,y,z]$ and $k[[x,y,z]]$. Similarly, 
we write $g_{4} = x^{2} + a(y,z)x + b(y,z)$, where $a(y,z) = -y^{3-r}z^{r+p-1}$, 
$b(y,z) = -y - y^{\frac{1}{2}(r-2)(3r-5)}z^{(-r^2+4r+p-2)}$. We can use 
Eisenstein's criterion with respect to the prime element $y$ to prove that 
$g_{4}$ is irreducible. 
The fact that $\{g_{4},g_{3}\}$ and $\{g_{3},g_{4}\}$ form regular sequences 
follow from the fact that $g_{3}$ and $g_{4}$ are both irreducible.\qed
\medskip

\begin{theorem} $P_{1}=I=\langle g_{1},g_{2},g_{3},g_{4}\rangle$. 
\end{theorem}

\proof We have $\rho_{1}(g_{i})=0$, for $1\leq i\leq 4$, hence $I\subset P_{1}$. 
Again by Proposition \ref{irreducible}, $\{g_{3},g_{4}\}\subset I$ forms a regular 
sequence of length $2$. Consider the chain 
$0\subset \langle g_{3}\rangle \subset I\subset P_{1}$. Given that, 
$\mathrm{ht}(P_{1})= 2$, $I\subset P_{1}$ and that $I$ is a prime ideal, 
we have $P_{1}=I$.\qed
\medskip

\begin{theorem}
\begin{enumerate}
\item If $r\in\{1,2\}$, then $I=\langle g_{3},g_{4}\rangle$.
\item if $r=3$, then $I=\langle g_{2}, g_{3},g_{4}\rangle$.
\end{enumerate}
\end{theorem}

\proof It is easy to see that $g_{1}= z(g_{3})-y(g_{2})$ and the proof 
of (2) follows easily. Let $r\in\{1,2\}$. We have 
\begin{align*}
g_{2}&=y(x^{2}-y-xy^{(3-r)}z^{(r-1+p)}-y^{\frac{1}{2}(r-2)(3r-5) }z^{(-r^2+4r-2+p)})-\\
&x(xy-z-y^{(4-r)}z^{(r-1+p)})
\end{align*}
Moreover, $\frac{3r^2-11r+12}{2}=3-r$ and $r^2-4r+2=-r$. 
Therefore $g_{2}=y(g_{4})-x(g_{3})$. Hence, in this case, $I=\langle g_{3},g_{4}\rangle$. \qed

\section{The ideal $\mathfrak{p}_{1}$ in the ring $k[[x,y,z]]$}
We now use Moh's techniques to 
prove that the set $\{g_{3},g_{4}\}$ forms a generating set for the ideal 
$\mathfrak{p}_{1}$ in the power series ring $k[[x,y,z]]$. 

\begin{theorem}
The set $\{g_{3},g_{4}\}$ forms a minimal generating set for the ideal $\mathfrak{p}_{1}$ in the ring $k[[x,y,z]]$
\end{theorem}

\proof We consider the case $n=1$, therefore we prove the followings,
\begin{enumerate}[(i)]
\item The $\sigma$-leading form of $g_{4}\in V_{2}$,
\item $x$( the $\sigma$-leading form of $g_{4}$) and the $\sigma$-leading form of 
$g_{3}$ generate $V_{3}$.
\end{enumerate}

\noindent\textbf{Case $r=1$.} We have
$$g_{3}=xy-(z+y^{3}z^{p}), \quad g_{4}=x^{2}-y-xy^{2}z^{p}-yz^{(1+p)}.$$
Therefore the $\sigma$-leading form of $g_{3}$ is $xy-z=h_{2}\in V_{3}$ and 
the $\sigma$-leading form of $g_{4}$ is $x^{2}-y=h_{1}\in V_{2}$. 
It is enough to show that $x(x^{2}-y)=x^{3}-xy$ and $xy-z$ generates $V_{3}$. 
Let $h(x,y,z)\in V_{3}$, then, 
$$h(x,y,z)\in W_{3}=\{\sigma\mbox{-homogeneous\, form\, of} \, \sigma\mbox{-order} \, 3\}$$
and it is also the $\sigma$-leading form of an element of $\mathfrak{p}_{1}$. As 
$h(x,y,z)\in W_{3}$, $h(x,y,z)$ is a $\sigma$-homogeneous form and 
$\mbox{ord}(\sigma(h(x,y,z)))=3$. Now any homogeneous total degree $3$ 
polynomials in $k[[x,y,z]]$ can be written as 
$$a_{1}x^{3}+a_{2}y^{3}+a_{3}z^{3}+a_{4}xy^{2}+a_{5}x^{2}y+a_{6}xyz+a_{7}xz^{2}+a_{8}x^{2}z+a_{9}yz^{2}+a_{10}y^{2}z.$$
We write the above expression as follows,
\begin{align*}
\sigma(h(x,y,z)) & = h(x,y^{2},z^{3})\\
&=a_{1}x^{3}+a_{2}y^{3}+a_{3}z^{3}+a_{4}xy^{2}+a_{5}x^{2}y \\
&+a_{6}xyz+a_{7}xz^{2}+a_{8}x^{2}z+a_{9}yz^{2}+a_{10}y^{2}z.
\end{align*}
We must have, $a_{2}=a_{5}=a_{6}=a_{7}=a_{8}=a_{9}=a_{10}=0$. Therefore, 
$h(x,y^{2},z^{3})=a_{1}x^{3}+a_{3}z^{3}+a_{4}xy^{2}$. Hence 
$h(x,y,z)=a_{1}x^{3}+a_{3}z+a_{4}xy$, since $h$ is $\sigma$-homogeneous. 
\medskip

Since $h(x,y,z)$ is a $\sigma$-leading form of an element of $\mathfrak{p}_{1}$, 
there exists $g(x,y,z)\in \mathfrak{p}_{1}$ such that the leading form of 
$\sigma(g(x,y,z))$ is $h(x,y,z)$. Let $g(x,y,z)=h(x,y,z)+f(x,y,z)$. Our claim is 
that, after applying $\rho_{1}$, any monomial in support of $f(x,y,z)$ cannot produce 
terms of the form $at^{3}$. If not, then we may write 
$\rho_{1}(a_{i_{1}i_{2}i_{3}}x^{i_{1}}y^{i_{2}}z^{i_{3}})=at^{3}+()$, and that forces 
$i_{1}+2i_{2}+3i_{3}=3$. The possible solution of this equation is 
$\{(1,1,0),(0,0,1), (3,0,0)\}$, hence the corresponding monomial would 
be $a_{(1,1,0)}xy$, $a_{(0,0,1)}z$, $a_{(3,0,0)}x^{3}$; which gives a contradiction as 
$h(x,y,z)$ is $\sigma$-homogeneous of order $3$ and the 
$\sigma$-leading form of $g(x,y,z)$. Thus our claim is proved.
\medskip

We apply $\rho_{1}$ on $ g(x,y,z)=h(x,y,z)+f(x,y,z)$ thus we get
$$0=a_{1}(t+t^{1+3(1+p)})^{3}+a_{3}t^{3}+a_{4}(t+t^{1+3(1+p)})t^{2}+\rho (f(x,y,z)).$$
Equating coefficient of $t^{3}$, we get, $a_{1}+a_{3}+a_{4}=0.$
We have,
\begin{align*}
h(x,y,z)&=a_{1}x^{3}+a_{3}z+a_{4}xy\\
&=a_{1}(x^{3}-xy)+a_{3}(z-xy)+(a_{1}+a_{3}+a_{4})xy\\
&=a_{1}(x^{3}-xy)+a_{3}(z-xy)
\end{align*}
Therefore $x^{3}-xy$ and $xy-z$ generate $V_{3}$.
\medskip

\noindent\textbf{Case $r=2$.} We have $$g_{3}=xy-[z+y^{2}z^{(1+p)}], \quad 
g_{4}=x^{2}-y-xyz^{(1+p)}-z^{(p+2)}.$$
The $\sigma$-leading forms of $g_{3}$ and $g_{4}$ are $(xy-z)$ and $(x^{2}-y)$ 
respectively; hence similar proof works.
\medskip

\noindent\textbf{Case $r=3$.} We have
$$g_{3}=xy-[z+yz^{(2+p)}], \quad g_{4}=x^{2}-y-xz^{(2+p)}-y^{2}z^{(1+p)}.$$
The $\sigma$-leading forms of $g_{3}$ and $g_{4}$ are $xy-z$ and $x^{2}-y$ respectively; 
hence similar proof works.
\medskip

Finally, to show that $\{g_{3},\, g_{4}\}$ forms minimal generating set for the ideal $\mathfrak{p}_{1}$ in the ring $k[[x,y,z]]$, we argue as follows. Height of the ideal $\mathfrak{p}_{1}$ is $2$ and the set $\{g_{3},\, g_{4}\}$ generates $\mathfrak{p}_{1}$. 
Hence by Krull's theorem the set $\{g_{3},\, g_{4}\}$ forms minimal generating set for the ideal $\mathfrak{p}_{1}$.\qed

\section{Explicit generators for $\mathfrak{p}_{n}$ in $k[[x,y,z]]$}
We now create a minimal generating set consisting of $n+1$ elements for the ideal $\mathfrak{p}_{n}$. 
We do that in two steps. First we write the 
$\sigma$-leading forms $f^{\sigma}$ and then we write the tails $f^{\tau}$ in order to create 
the entire polynomial $f = f^{\sigma} + f^{\tau}$. While writing these, we will first write the 
expressions as linear combinations of suitable 
monomials and then calculate the coefficients through a system of linear equations. 
Moh's theorem 
will be used for this purpose, which has been described in Section 1. 
We indicate the main steps of the procedure below, to be explained in the rest of the paper.
\medskip

\begin{enumerate}
\item \textbf{Step 1.} Write $f_{1}^{\sigma}$ and $f_{2}^{\sigma}$.
\medskip

\item \textbf{Step 2.} Form the equation \ref{geneqn} and write the systems of 
equations \ref{geneqn8} to \ref{geneqn12}, if $i$ is odd, and equations \ref{geneqn14} 
to \ref{geneqn18}, if $i$ is even. 
\medskip

\item \textbf{Step 3.} Compute the coefficients $a_{(u,v)}$ and $c_{(p,q)}$ from the 
above systems of equations. Existence of nontrivial solutions for $c_{(p,q)}$ is 
guaranteed. 
\medskip

\item \textbf{Step 4.} Form the equation \ref{geneqn19} with the same coefficients 
$a_{(p,q)}$, form the systems of equations based on odd and even values of $i$.
\medskip

\item \textbf{Step 5.} Compute the coefficients $d_{(r,s)}$ in a similar fashion from 
the said systems.
\medskip

\item \textbf{Step 6.} Write  $f_{i}=f_{i}^{\sigma}+f_{i}^{\tau}$, $1\leq i\leq n+1$.
\end{enumerate}
\medskip

\begin{center}
\textbf{The $\mathbf{\sigma}$ - leading forms} 
\end{center}
Computations with some special cases help us guess finitely many polynomials, 
which are the expected $\sigma$-leading forms of some prospective generators 
of $\mathfrak{p}_{n}$. For $c(i,j)\in k$, we write,
\begin{align*}
f_{2i-1}^{\sigma}&= \displaystyle\sum_{j=1}^{i}(-1)^{(j-1)} c_{(2i-1,j)}x^{(2m-i-j+2)}y^{(2j-2)}z^{(i-j)}\\
& + \displaystyle\sum_{j=1}^{m+1-i}(-1)^{(j+i-1)} c_{(2i-1,i+j)} x^{(m+1-i-j)}y^{(2j-1)}z^{(m-1+i-j)},\quad 1\leq i\leq m
\end{align*}
\begin{align*}
 f_{2i}^{\sigma}&= \displaystyle\sum_{j=1}^{i} (-1)^{(j-1)} c_{(2i,j)} x^{(2m-i-j+1)}y^{(2j-1)}z^{(i-j)}\\
 & + \displaystyle\sum_{j=1}^{m+1-i} (-1)^{(j+i-1)} c_{(2i,i+j)} x^{(m+1-i-j)}y^{(2j-2)}z^{(m+i-j)},\quad 1\leq i\leq m. 
\end{align*}
Let us assume that, $c_{(1,1)}=c_{(2,1)}=1$ and $c_{(1,j+1)}=\binom{m+j-1}{j-1}\binom{2m}{m-j}$, $c_{(2,j+1)}=\binom{m+j-2}{j-1}\binom{2m-1}{m-j}$ for $1\leq j\leq m$. Therefore, 
\begin{align*}
&f_{1}^{\sigma}= x^{2m}+\displaystyle\sum_{j=1}^{m} (-1)^{j}\binom{m+j-1}{j-1}\binom{2m}{m-j} x^{(m-j)}y^{(2j-1)}z^{(m-j)}\\
&f_{2}^{\sigma}=  x^{2m-1}y+\displaystyle\sum_{j=1}^{m}(-1)^{j}\binom{m+j-2}{j-1}\binom{2m-1}{m-j} x^{(m-j)}y^{(2j-2)}z^{(m+1-j)}
\end{align*}
\medskip

\noindent Our aim is to show that there exist polynomials 
$f_{1},\ldots,f_{n+1}\in \mathfrak{p}_{n}$ with $\sigma$-leading forms 
$f_{1}^{\sigma},\ldots,f_{n+1}^{\sigma}$, and that 
$f_{1},\ldots,f_{n+1}$ form a minimal set of generators of $\mathfrak{p}_{n}$. 
In order to calculate the coefficients $c_{(p,q)}$ of $f_{i}^{\sigma}$ for $3\leq i\leq n+1$, 
Moh's work suggests us to consider the following equations,
\begin{align}\label{geneqn}
a_{(i,i)}zf_{i}^{\sigma}+a_{(i,(i+1))}yf_{i+1}^{\sigma}+a_{(i,(i+2))}xf_{i+2}^{\sigma}=0
\end{align}
for $1\leq i\leq n-1$ and $a_{(i,i)},a_{(i,(i+1))},a_{(i,(i+2))}\in k$ with $a_{(i,(i+2))}\neq 0$. 
\medskip

First we consider equation \ref{geneqn} for odd values of $i$, where we write $2i-1$ instead 
of introducing a new index and obtain the equation,
\begin{align}\label{geneqn7}
a_{(2i-1,2i-1)}zf_{2i-1}^{\sigma}+a_{(2i-1,2i)}yf_{2i}^{\sigma}+a_{(2i-1,2i+1)}xf_{2i+1}^{\sigma}=0.
\end{align} 

\noindent Equating coefficients of monomials in equation \ref{geneqn7} we obtain the following system of equations \ref{geneqn8} to \ref{geneqn12}:
\begin{align}\label{geneqn8}
 a_{(2i-1,2i-1)}c_{(2i-1,1)}+a_{(2i-1,2i+1)}c_{(2i+1,1)}=0,
\end{align}
\begin{align}\label{geneqn9}
a_{(2i-1,2i-1)}c_{(2i-1,j)}-a_{(2i-1,2i)}c_{(2i,j-1)}+a_{(2i-1,2i+1)}c_{(2i+1,j)}  =0; \, 2\leq j\leq i;
\end{align}
\begin{align}\label{geneqn10}
a_{(2i-1,2i)}c_{(2i,i)}-a_{(2i-1,2i+1)}c_{(2i+1,i+1)}=0
\end{align}
\begin{align}\label{geneqn11}
a_{(2i-1,2i-1)}c_{(2i-1,i+j)}+a_{(2i-1,2i)}c_{(2i,i+j)}-a_{(2i-1,2i+1)}c_{(2i+1,i+j+1)}=0; 
\, 1\leq j\leq m-i;
\end{align}
\begin{align}\label{geneqn12}
a_{(2i-1,2i-1)}c_{(2i-1,m+1)}+a_{(2i-1,2i)}c_{(2i,m+1)}=0.
\end{align}
\medskip

Next we consider equation \ref{geneqn} for even values of $i$, where we write $2i$ instead 
of introducing a new index and obtain,
\begin{align}\label{geneqn13}
a_{(2i,2i)}zf_{2i}^{\sigma}+a_{(2i,2i+1)}yf_{2i+1}^{\sigma}+a_{(2i,2i+2)}xf_{2i+2}^{\sigma}=0.
\end{align}

\noindent Equating coefficients of monomials in equation \ref{geneqn13} we get the following system of equations \ref{geneqn14} to \ref{geneqn18}:
\begin{align}\label{geneqn14}
 a_{(2i,2i)}c_{(2i,j)}+a_{(2i,2i+1)}c_{(2i+1,j)}+a_{(2i,2i+2)}c_{(2i+2,j)}=0;\, 1\leq j\leq i
\end{align}
\begin{align}\label{geneqn15}
 a_{(2i,2i+1)}c_{(2i+1,i+1)}+a_{(2i,2i+2)}c_{(2i+2,i+1)}=0;
\end{align}
\begin{align}\label{geneqn16}
 a_{(2i,2i)}c_{(2i,i+1)}-a_{(2i,2i+2)}c_{(2i+2,i+2)}=0;
\end{align}
\begin{align}\label{geneqn17}
 a_{(2i,2i)}c_{(2i,i+j)}+a_{(2i,2i+1)}c_{(2i+1,i+j)}-a_{(2i,2i+2)}c_{(2i+2,i+j+1)}=0, \, 2\leq j\leq m-i;
\end{align}
\begin{align}\label{geneqn18}
 a_{(2i,2i)}c_{(2i,m+1)}+a_{(2i,2i+1)}c_{(2i+1,m+1)}=0.
\end{align} 
\medskip

\noindent\textbf{The cases $i=1$ and $i=2$:} Let us take $i=1$. The equation \ref{geneqn} takes the form,
\begin{align}\label{en1}
a_{(1,1)}zf_{1}^{\sigma}+a_{(1,2)}yf_{2}^{\sigma}+a_{(1,3)}xf_{3}^{\sigma}=0.
\end{align}
We know $f_{1}^{\sigma}$ and $f_{2}^{\sigma}$ and our aim is to compute $f_{3}^{\sigma}$. 
Equating coefficients of monomials we get, the system of equations \ref{en2} to \ref{en5}:
\begin{align}\label{en2}
a_{(1,1)}+a_{(1,3)}c_{(3,1)}=0;
\end{align}
\begin{align}\label{en3}
a_{(1,2)}-a_{(1,3)}c_{(3,2)}=0;
\end{align}
\begin{align}\label{en4}
a_{(1,1)}\binom{m+j-1}{j-1}\binom{2m}{m-j}-a_{(1,2)}\binom{m+j-2}{j-1}\binom{2m-1}{m-j}-a_{(1,3)}c_{(3,j+2)}=0; \, 1\leq j\leq m-1;
\end{align}
\begin{align}\label{en5}
a_{(1,1)}\binom{2m-1}{m-1}+a_{(1,2)}\binom{2m-2}{m-1}=0.
\end{align}

\noindent For any non-zero value for $a_{(1,3)}$, say $0\neq a=a_{(1,3)}$, we have 
$a_{(1,1)}=-a c_{(3,1)}$ and $a_{(1,2)}=ac_{(3,2)}$ from equations \ref{en2} and \ref{en3}. 
Substituting the values of $a_{(1,1)},a_{(1,2)},a_{(1,3)}$ in equations \ref{en4}, 
\ref{en5} we get $m$ equations in $m+1$ unknowns $c_{(3,j)},\, 1\leq j\leq m+1$. 
Therefore nontrivial solutions exists. We pick a non trivial solution, which 
gives us the desired coefficients of $f^{\sigma}_{3}$. Next we consider $i=2$, that is, 
\begin{align}\label{en1}
a_{(2,2)}zf_{2}^{\sigma}+a_{(2,3)}yf_{3}^{\sigma}+a_{(2,4)}xf_{4}^{\sigma}=0,\quad a_{(2,4)}\neq 0.
\end{align} 
Having known the coefficients of $f_{2}^{\sigma}$ and $f_{3}^{\sigma}$, we can calculate 
all the coefficients of $f_{4}^{\sigma}$ from the system of equations in a similar 
fashion. Existence of nontrivial solution is once again guaranteed because 
there are $m$ equations in $m+1$ unknowns. 
\medskip

\noindent\textbf{The case for arbitrary $i$:} Let us now return to the general case, 
where we would be using an inductive argument with the starting step being $i=1$. 
Let us present the argument only for odd values of $i$; the same for even values 
of $i$ would be similar because of similarity in the system of equations. We may 
assume that $c_{i,1}=1$ for $1\leq i\leq n+1$. Then, from equation \ref{geneqn8} 
we get, 
$$a_{(2i-1,2i-1)} + a_{(2i-1,2i+1)} = 0.$$ 
We take a non zero value for $a_{(2i-1,2i+1)}$, say $a_{(2i-1,2i+1)}=1$; then 
$a_{(2i-1,2i-1)}=-1$. If $c_{(2i,m+1)}=0$, then we may take any non zero value 
for $a_{(2i-1,2i)}$ and if $c_{(2i,m+1)}\neq 0$ then we take 
$a_{(2i-1,2i)}=\dfrac{c_{(2i-1,m+1)}}{c_{(2i,m+1)}}$. We substitute these 
values in the equations \ref{geneqn9}, \ref{geneqn10} and \ref{geneqn11} 
and get a solution for $c_{(2i+1,j)}$, for $2\leq j\leq m+1$. Note that $c_{(2i-1,j)}$ and 
$c_{(2i,j)}$ are known from the previous step. Thus, we can find coefficients of 
$f_{i}^{\sigma}$ using the relation \ref{geneqn} and corresponding equations. 
This therefore gives us the desired $\sigma$-leading forms of the prospective generators 
explicitly. We have illustrated this and the subsequent steps Example \ref{mainex}. 
\bigskip

\begin{center}
\textbf{The tail $f^{\tau}$}
\end{center}  
We now proceed to write the tails $f_{1}^{\tau}, \ldots f_{n+1}^{\tau}$. 
Before we begin, a few technical lemmas are in order. Let us recall that, $\lambda>n(n+1)m$ and 
$\gcd(\lambda, m)=1$. Let us write, $\lambda=(nm+w)(n+1)+r$, then $1\leq r\leq n+1$, such that 
$\gcd(\lambda, m)=1$. Further, if we write $w=\alpha(n+2)+q$ and $w+n+1=\alpha^{'}(n+2)+q^{'}$, 
then $0\leq q$, $q^{'}\leq n+2$. 

\begin{lemma}
Suppose $t_{1} =(q+(r-1)(n+1)) \, \hbox{mod} (n+2)$ and $t_{1}^{'}=(q^{'}+(r-1)(n+1))\mathrm{mod}(n+2)$. 
Then \, $n+2$\, divides \, $n(m+1)+\lambda-(n+1)t_{1}+1$\, and \, $nm+\lambda-(n+1)t_{1}^{'}$.
\end{lemma}

\proof Let $q+(r-1)(n+1)-t_{1}=(n+2)k$. Substituting the values of 
$\lambda$ and $w$ we get,
\begin{align*}
n(m+1)+\lambda-(n+1)t_{1}+1&=\frac{n^2+3n+2}{2}-(n+1)t_{1}+\lambda \\
&=(n+2)[m-2mk-rn+2\alpha m+2(m-1)^2].
\end{align*}
Let $q^{'}+(r-1)(n+1)-t^{'}_{1}=(n+2)k^{'}$; substituting the values of 
$\lambda$ and $w$ we get: $\, nm+\lambda-(n+1)t_{1}^{'}=(n+2)[mn-k^{'}-rn+2m\alpha^{'}]$. \qed 
\medskip

\begin{lemma}\label{binomial}
Let 
\begin{enumerate}
\item[(i)] $\nu_{l}=\binom{2m}{m+l}-\binom{m-1}{0}\binom{m}{l}-\binom{m}{1}\binom{m-1}{l}-\binom{m+1}{2}\binom{m-2}{l}-\cdots-\binom{2m-(l+1)}{m-l}\binom{m-(m-l)}{l}$, where $0\leq l \leq m$;
\medskip

\item[(ii)] $\mu_{l}=\binom{2m}{l}-\binom{m}{0}\binom{2m}{m-1}\binom{m-1}{l}+\binom{m+1}{1}\binom{2m}{m-2}\binom{m-2}{l}-\binom{m+2}{2}\binom{2m}{m-3}\binom{m-3}{l}+\cdots+(-1)^{(m-l)}\binom{2m-1-l}{m-1-l}\binom{2m}{l}\binom{m-(m-l)}{l}$, 
where $0\leq l \leq m-1$.
\end{enumerate}

Then $\mu_{l}=0$ and $\nu_{l}=0$, for all values of $l$. 
\end{lemma}

\proof (i). We use the identity $$ \displaystyle\sum_{k=0}^{t}\binom{a+k}{k}\binom{s+t-k}{t-k}=\binom{a+s+t+1}{t},$$  
and substitute $s=l, t=m-l, a=m-1$.
\medskip

\noindent (ii). Using the identity $\binom{r}{s}\binom{s}{t}=\binom{r}{t}\binom{r-t}{s-t}$, we rewrite 
$\mu_{l}$ as 
$$\mu_{l} = \binom{2m}{l}-\binom{2m}{l}\cdot \left[\displaystyle{\sum_{k=1}^{m-l}}(-1)^{k-1}\binom{2m-l}{m+k}\binom{m+k-1}{k-1}\right].$$ 
Therefore, it is enough to prove that 
$$\displaystyle{\sum_{k=1}^{m-l}}(-1)^{k-1}\binom{2m-l}{m+k}\binom{m+k-1}{k-1}=1.$$ 
We prove this by induction on $l$. For $m-l=1$, we get $\binom{m+1}{m+1}\binom{m}{0}=1$. 
Let the result hold for $m-l=u$, that is, 
$$\displaystyle{\sum_{k=1}^{u}}(-1)^{k-1}\binom{2m-l}{m+k}\binom{m+k-1}{k-1}=\displaystyle{\sum_{k=1}^{u}}(-1)^{k-1}\binom{m+u}{m+k}\binom{m+k-1}{k-1}=1.$$
For $m-l=u+1$, using the result $\binom{r}{s}=\binom{r-1}{s-1}+\binom{r-1}{s}$, we get\\
\begin{align*}
\displaystyle{\sum_{k=1}^{u+1}}(-1)^{k-1}\binom{2m-l}{m+k}\binom{m+k-1}{k-1} = \displaystyle{\sum_{k=1}^{u+1}}(-1)^{k-1}\binom{m+u+1}{m+k}\binom{m+k-1}{k-1}\\
= \displaystyle{\sum_{k=1}^{u}}(-1)^{k-1}\binom{m+u}{m+k}\binom{m+k-1}{k-1}
+ \displaystyle{\sum_{k=1}^{u+1}}(-1)^{k-1}\binom{m+u}{m+k-1}\binom{m+k-1}{k-1}
\end{align*}
We have,  
\begin{align*}
\displaystyle{\sum_{k=1}^{u+1}}(-1)^{k-1}\binom{m+u}{m+k-1}\binom{m+k-1}{k-1}\\
= \binom{m+u}{m}\displaystyle{\sum_{k=1}^{u+1}}(-1)^{k-1}\binom{u}{k-1} =0.
\end{align*}
This proves the desired result. \qed
\medskip

Let us write $t_{l}=(t_{l-1}+2)\mathrm{mod}(n+2)$ and 
$t^{'}_{l}=(t^{'}_{l-1}+2)\mathrm{mod}(n+2)$, $l\geq 2$. We now 
define a bunch of monomials which would appear in the expression of the tails of the 
prospective generators. Given $1\leq l\leq m$, we define
\begin{itemize}
\item $m_{((n-2l+3),1)}=x^{l-1}y^{t_{1}}z^{(\frac{n(m+1)+\lambda-(n+1)t_{1}+1}{n+2}+m-l)}$.
\medskip

\item Having defined the monomials $m_{((n-2l+3),s-1)}$,  $2\leq s \leq l$, we define  
\begin{align*}m_{((n-2l+3),s))} & = 
\begin{cases} x^{l-s}y^{t_{s}}z^{\gamma_{((n-2l+3),s-1)}+t_{s-1}-t_{s}} & \mbox{if} \quad t_{s-1}+2\geq n+2,\\
x^{l-s}y^{t_{s}}z^{\gamma_{((n-2l+3),s-1)}-1} & \mbox{if} \quad t_{s-1}+2< n+2;
\end{cases}
\end{align*}
where $\gamma_{((n-2l+3),s-1)}$ is the power of $z$ in $m_{((n-2l+3),s-1)}$.
\medskip

\item $m_{((n-2l+2),1)}=x^{l}y^{t^{'}_{1}}z^{(\frac{nm+\lambda-(n+1)t_{1}^{'}}{n+2}+m-l)}$.
\medskip

\item Having defined the monomials $m_{((n-2l+2),s-1)}$ for $2\leq s\leq l+1$, we define  
\begin{align*}
m_{((n-2l+2),s)} & = 
\begin{cases}
x^{(l-s+1)}y^{t^{'}_{s}}z^{(\gamma_{((n-2l+2),s-1)}^{'}+t_{s-1}^{'}-t^{'}_{s})} & 
\mbox{if} \quad t^{'}_{s-1}+2\geq n+2,\\
x^{(l-s+1)}y^{t^{'}_{s}}z^{(\gamma^{'}_{((n-2l+2),s-1)}-1)} & \mbox{if} \quad t^{'}_{s-1}+2< n+2;
\end{cases}
\end{align*}
where $\gamma^{'}_{((n-2l+2),s-1)}$ is the power of $z$ in $m_{((n-2l+2),s-1)}$.
\end{itemize}
\medskip

\noindent We now have to prove that the integers $\gamma_{((n-2l+3),s-1)}$ and 
$\gamma^{'}_{((n-2l+2),s-1)}$, occurring as powers of $z$ in the monomials 
defined above are indeed integers. Since $\lambda>nm(n+1)$ and 
$0\leq t_{1},t_{1}^{'}\leq n-1$, we have 
$$\frac{n(m+1)+\lambda-(n+1)t_{1}+1}{n+2}>m+1$$ and 
$$\frac{nm+\lambda-(n+1)t_{1}^{'}}{n+2}>m+1.$$ Therefore $\gamma_{((n-2l+3),s)}>0$ 
and $\gamma^{'}_{((n-2l+2),s)}>0 $ for all values of $l,s$.
\medskip

We are now ready to define the tails. Let us first write $f_{1}^{\tau}$ and $f_{2}^{\tau}$ 
as $f_{1}^{\tau}=-\displaystyle\sum_{s=1}^{m+1}\binom{m+s-2}{s-1} m_{(1,s)}$ 
and $f_{2}^{\tau}=-\displaystyle\sum_{s=1}^{m}\binom{m+s-2}{s-1} m_{(2,s)}$. 
The polynomials $f_{1}$ and $f_{2}$ are defined as 
\begin{align*}
f_{1}=f_{1}^{\sigma}+f_{1}^{\tau}=& x^{2m}+\displaystyle\sum_{j=1}^{m} (-1)^{j}\binom{m+j-1}{j-1}\binom{2m}{m-j} x^{(m-j)}y^{(2j-1)}z^{(m-j)}\\
&-\displaystyle\sum_{s=1}^{m+1}\binom{m+s-2}{s-1} m_{(1,s)};
\end{align*} 
\begin{align*}
 f_{2}=f_{2}^{\sigma}+f_{2}^{\tau}= & x^{2m-1}y+\displaystyle\sum_{j=1}^{m} (-1)^{j}\binom{m+j-2}{j-1}\binom{2m-1}{m-j} x^{(m-j)}y^{(2j-2)}z^{(m+1-j)}\\
 &-\displaystyle\sum_{s=1}^{m}\binom{m+s-2}{s-1} m_{(2,s)}.
\end{align*}  
\medskip

\noindent For $1\leq l\leq m $, let us write 
$$f_{(n-2l+3)}^{\tau}=-\left(d_{((n-2l+3),1)}m_{((n-2l+3),1)}+\cdots +d_{((n-2l+3),l)} m_{((n-2l+3),l)}\right),\quad $$
and  $$f_{n-2l+2}^{\tau}=-\left(d_{((n-2l+2),1)}m_{(n-2l+2,1)}+\cdots +d_{((n-2l+2),(l+1))} m_{((n-2l+2),(l+1))}\right).$$
Note that so far we have only written the monomials appearing in the tails 
$f_{3}^{\tau}, \ldots, f_{n+1}^{\tau}$. We now have to compute the coefficients 
attached to each monomial and we follow a similar strategy like determining the 
coefficients $d_{(r,s)}$ of $f_{i}^{\tau}$ inductively through the relation 
\begin{align}\label{geneqn19}
a_{(i,i)}zf_{i}^{\tau}+a_{(i,(i+1))}yf_{i+1}^{\tau}+a_{(i,(i+2))}xf_{i+2}^{\tau}=0, \quad \mbox{for} \quad 
1\leq i\leq n-1,
\end{align}
where $a_{(i,i)},a_{(i,i+1)},a_{(i,i+2)}\in k$ are the same as in the equation \ref{geneqn}. 
We write $f_{i}=f_{i}^{\sigma}+f_{i}^{\tau}$, $1\leq i\leq n+1$, and it follows from 
equations \ref{geneqn} and \ref{geneqn19} that 
\begin{align}\label{geneqn*}
a_{(i,i)}zf_{i}+a_{(i,(i+1))}yf_{i+1}+a_{(i,(i+2))}xf_{i+2}=0, \quad \mbox{for} \quad 
1\leq i\leq n-1.
\end{align}
The most important observation after defining $f_{1}, \ldots , f_{n+1}$ is presented 
in the following Lemma and the proof follows easily from the definition of $\sigma$-weight.
\medskip

\begin{lemma}$f_{1}^{\sigma},\ldots,f_{n+1}^{\sigma}$ are indeed the $\sigma$-leading forms of $f_{1},\ldots,f_{n+1}$.
\end{lemma}
\medskip

\begin{example}\label{mainex} We now understand the main steps involved in the construction 
illustrated in the example below. One has to understand that the construction of the 
generating polynomials $f_{1}, \ldots , f_{n+1}$ is purely based on the equations 
\ref{geneqn} and \ref{geneqn19}. Finally, equation \ref{geneqn*} will be used along with 
Moh's theorem to show that $f_{1}, \ldots , f_{n+1}$ indeed generate $\mathfrak{p}_{n}$. 
\medskip

Let $n=3$, $m=(n+1)/2=2$, $\lambda = 27$. 
Evidently, $\lambda > 24$ and $\gcd (\lambda, m)=1$. Moh's curve is 
defined by the parametrization $\rho_{3}(x)=t^{6}+t^{33}$, 
$\rho_{3}(y)=t^{8}$, $\rho_{3}(z)=t^{10}$. We show how to write 
a generating set $\{f_{1}, f_{2}, f_{3}, f_{4}\}$ of $\mathfrak{p}_{3}$, 
starting with explicit expressions for the $\sigma$- homogeneous 
terms $f_{i}^{\sigma}$, $i=1, 2, 3, 4$ followed by the tails 
$f_{i}^{\tau}$, $i=1, 2, 3, 4$. We start with
\begin{itemize}
\item $f_{1}^{\sigma}=x^{4}-4xyz+3y^{3}$,
\item $f_{2}^{\sigma}=x^{3}y-3xz^{2}+2y^{2}z$,
\end{itemize}
and then write
\begin{itemize}
\item $f_{3}^{\sigma}=c_{(3,1)}x^{3}z-c_{(3,2)}x^{2}y^{2}+c_{(3,3)}yz^{2}$,
\item $f_{4}^{\sigma}=c_{(4,1)}x^{2}yz-c_{(4,2)}xy^{3}+c_{(4,3)} z^{3}$.
\end{itemize}
We use equation \ref{geneqn} for $i=1$ and obtain,
$$a_{(1,1)}zf_{1}^{\sigma}+a_{(1,2)}yf_{2}^{\sigma}+a_{(1,3)}x f_{3}^{\sigma}=0.$$ 
Equating coefficients of the monomials we get the following system of equations: 
\begin{enumerate}
\item $a_{(1,1)}+a_{(1,3)}c_{(3,1)}=0$;
\item $-4a_{(1,1)}-3a_{(1,2)}+c_{(3,3)}a_{(1,3)}=0 $
\item $3a_{(1,1)}+2a_{(1,2)}=0$
\item $-a_{(1,3)}c_{(3,2)}+a_{(1,2)}=0$
\end{enumerate}
Let $a=a_{(1,3)}\neq 0$, we get $a_{(1,1)}=-ac_{(3,1)}$ and $a_{(1,2)}=-\frac{3}{2}ac_{(3,1)}$, also we have  $-c_{(3,1)}+2c_{(3,3)}=0$ and $-2c_{(3,2)}+3c_{(3,1)}=0$. 
Suppose we choose $c_{(3,1)}=2, c_{(3,2)}=3,c_{(3,3)}=1$, then $a_{(1,1)}=-2a$ and $a_{(1,2)}=3a$. Furthermore we choose $a=1$ then $a_{(1,1)}=-2, a_{(1,2)}=3,a_{(1,3)}=1$.
\medskip

\noindent Similarly using the equation \ref{geneqn} for $i=2$ we get, 
$$a_{(2,2)}zf_{2}^{\sigma}+a_{(2,3)}yf_{3}^{\sigma}+a_{(2,4)}x f_{4}^{\sigma}=0$$
from the above equation we get the following equations
\begin{enumerate}
\item $a_{(2,2)}+2a_{(2,3)}+a_{(2,4)}c_{(4,1)}=0$,
\item $3a_{(2,2)}-a_{(2,4)}c_{(4,3)}=0$,
\item $2a_{(2,2)}+a_{(2,3)}=0$,
\item $3a_{(2,3)}+a_{(2,4)}c_{(4,2)}=0$.
\end{enumerate}
Again we choose $a_{24}=1$ and proceed by same way we can take a solution,
$c_{(4,1)}=1, c_{(4,2)}=2, c_{(4,3)}=1$. Then $a_{(2,2)}=\frac{1}{3},a_{(2,3)}=-\frac{2}{3},a_{(2,4)}=1$.
\medskip

Now we write  $f_{i}^{\tau},\, 1\leq i\leq 4$. We write 
$\lambda=(nm+w)(n+1)+r$, and we get $w=0$ and $r=3$. Furthermore, 
if we write $w=\alpha(n+2)+q$, we get $\alpha=0$, $q=0$, $w+n+1=\alpha^{'}(n+2)+q^{'}$ and 
$q^{'}=4$, $\alpha^{'}=0$. Now using the formulae we calculate $t_{1}=3$ and $t_{1}^{'}=2$. 
We know that 
$$f_{1}^{\tau}=-\displaystyle\sum_{s=1}^{3}\binom{m+s-2}{s-1} m_{(1,s)};$$
and 
$$f_{2}^{\tau}=-\displaystyle\sum_{s=1}^{2}\binom{m+s-2}{s-1} m_{(2,s)}.$$

Let us take $l=2$ in $m_{((n-2l+2),1)}$ and substitute the value of $\lambda$, $t_{1}^{'}$; 
we get $m_{(1,1)}=x^{2}y^{2}z^{5}$. Taking $s=2$ in $m_{((n-2l+2),s)}$, we get the 
monomial $m_{(1,2)}$. We see that $t_{1}^{'}+2< 5$, therefore using the corresponding formulae and 
$t_{2}^{'}=4$, $\gamma_{(1,1)}=5$ we get $ m_{(1,2)}= xy^{4}z^{4}$. Similar calculation for 
$s=3$ yields $m_{(1,3)}= yz^{7}$. Hence, $f_{1}^{\tau}=-x^{2}y^{2}z^{5}-2xy^{4}z^{4}-3yz^{7}$. 
Let us now take $l=2$ in $m_{((n-2l+3),1)}$ and substitute the values of $\lambda$, $t_{1}$; 
we get $m_{(2,1)}=xy^{3}z^{5}$. Taking $s=2$ in $m_{((n-2l+3),s)}$, we get the monomial $m_{(2,2)}$. 
We see that $t_{1}+2=5$, therefore using the corresponding formulae and $t_{2}=0$, $\gamma_{(2,1)}=5$, 
we obtain $m_{(2,2)}= z^{8}$. Similarly calculation for $s=3$ yields $f_{2}^{\tau}=-xy^{3}z^{5}-2z^{8}$. 
Putting $l=1$ in $m_{((n-2l+2),1)}$ and $m_{((n-2l+3),1)}$ we get the monomials $m_{(3,1)}=xy^{2}z^{6}$ 
and $m_{(4,1)}=y^{3}z^{6}$. Now using the formulae for $m_{((n-2l+2),s)}$ we get $m_{(3,2)}=y^{4}z^{5}$. 
Hence $f_{3}^{\tau}=-d_{(3,1)}xy^{2}z^{6}-d_{(3,2)}y^{4}z^{5}$ and $f_{4}^{\tau}=-d_{(4,1)}y^{3}z^{6}$.
\medskip

Now we have 
\begin{itemize}
\item $f_{1}=f_{1}^{\sigma}+f_{1}^{\tau}=x^{4}-4xyz+3y^{3}-x^{2}y^{2}z^{5}-2xy^{4}z^{4}-3yz^{7}$,
\item $f_{2}=f_{2}^{\sigma}+f_{2}^{\tau}=x^{3}y-3xz^{2}+2y^{2}z-xy^{3}z^{5}-2z^{8}$,
\item $f_{3}=f_{3}^{\sigma}+f_{3}^{\tau}=2x^{3}z-3x^{2}y^{2}+yz^{2}-d_{(3,1)}xy^{2}z^{6}-d_{(3,2)}y^{4}z^{5}$,
\item $f_{4}=f_{4}^{\sigma}+f_{4}^{\tau}=x^{2}yz-2xy^{3}+ z^{3}-d_{(4,1)}y^{3}z^{6}$.
\end{itemize}

Using the equation $a_{(1,1)}zf_{1}+a_{(1,2)}yf_{2}+a_{(1,3)}x f_{3}=0$, we get $d_{(3,1)}=2,\, d_{(3,2)}=1$. 
Similarly using the equation $a_{(2,2)}zf_{2}+a_{(2,3)}yf_{3}+a_{(2,4)}x f_{4}=0$, we get $d_{(4,1)}=1$. 
Therefore,
$$f_{3}=f_{3}^{\sigma}+f_{3}^{\tau}=2x^{3}z-3x^{2}y^{2}+yz^{2}-2xy^{2}z^{6}-y^{4}z^{5},$$
$$f_{4}=f_{4}^{\sigma}+f_{4}^{\tau}=x^{2}yz-2xy^{3}+ z^{3}-y^{3}z^{6}.$$
\end{example}
\medskip

\begin{proposition}\label{f1f2}
The polynomials $f_{1},f_{2}$ belong to the kernel $\mathfrak{p}_{n}$.
\end{proposition}
\proof At first we show that $f_{1}$ lies in the kernel $\mathfrak{p}_{n}$. We claim that, 
$$\rho_{n}(m_{(1,u)})=\displaystyle\sum_{j=0}^{m+1-u} \binom{m+1-u}{j}t^{nm(n+1)+\lambda(m+j)}.$$ We proceed by induction on $u$. For $u=1$, it straightforward that, $$\rho_{n}(m_{(1,1)})=\displaystyle\sum_{j=0}^{m} \binom{m}{j}t^{nm(n+1)+\lambda(m+j)}.$$
Let us assume by induction it is true for $u=s$, then we have relation, 
 $$mn(m+1-s)+\lambda j +m(n+2)[\gamma_{(1,s-1)}^{'}+t^{'}_{s-1}]-mt^{'}_{s}=mn(n+1)+\lambda(m+j).$$
 We have to show that,  $$mn(m-s)+\lambda j +m(n+2)[\gamma_{(1,s)}^{'}+t^{'}_{s}]-mt^{'}_{s+1}=mn(n+1)+\lambda(m+j).$$ Where, $\gamma_{(1,s)}^{'}=\gamma_{(1,s-1)}^{'}+(t^{'}_{s-1}-t^{'}_{s} )$ and $(t^{'}_{s}-t^{'}_{s+1})=n$.
 We have,
 \begin{align*}
 & mn(m-s)+\lambda j +m(n+2)[\gamma_{(1,s)}^{'}+t^{'}_{s}]-mt^{'}_{s+1}\\
 &=mn(m-s)+\lambda j +m(n+2)[\gamma_{(1,s-1)}^{'}+t^{'}_{s-1}-t^{'}_{s}+t^{'}_{s}]-mt^{'}_{s+1}\\
 &=mn(m-s+1)+\lambda j +m(n+2)[\gamma_{(1,s-1)}^{'}+t^{'}_{s-1}]-mt^{'}_{s}+m(t^{'}_{s}-t^{'}_{s+1})-mn\\
 &=mn(n+1)+\lambda(m+j).
 \end{align*}
 Thus our claim is proved. Next we have,
\begin{align*}
\rho_{n}(f_{1})=& (t^{nm}+t^{(nm+\lambda)})^{2m}\\
&+\displaystyle\sum_{j=1}^{m} (-1)^{j}\binom{m+j-1}{j-1}\binom{2m}{m-j} (t^{nm}+t^{(nm+\lambda)})^{(m-j)}t^{mn(m+j)}\\
&-\displaystyle\sum_{s=1}^{m+1}\binom{m+s-2}{s-1} \rho_{n}(m_{(1,s)}).
\end{align*} 
Therefore,
\begin{align*}
\rho_{n}(f_{1})=& (t^{nm}+t^{(nm+\lambda)})^{2m}\\
&+\displaystyle\sum_{j=1}^{m} (-1)^{j}\sum_{v=0}^{m-j}\binom{m+j-1}{j-1}\binom{2m}{m-j}\binom{m-j}{v} t^{(nm(m-j)+\lambda v)}t^{mn(m+j)}\\
&-\displaystyle\sum_{s=1}^{m+1}\binom{m+s-2}{s-1} \rho_{n}(m_{(1,s)}).
\end{align*} 
Hence $$\rho_{n}(f_{1})=\displaystyle\sum_{j=0}^{m-1}\mu_{j}t^{(nm(n+1)+\lambda j)}+\sum_{j=0}^{m}\nu_{j}t^{(nm(n+1)+\lambda (m+j))}=0.$$
Similarly we can show that, $f_{2}$ also lies in $\mathfrak{p}_{n}$.\qed
\medskip

\begin{theorem}
The polynomials $f_{1},\ldots,f_{n+1}$ defined above generate $\mathfrak{p}_{n}$ minimally.
\end{theorem}

\proof Moh proved that $\mathfrak{p}_{n}$ is minimally generated by $n+1$ 
polynomials. Therefore, it is enough to prove that $f_{1},\ldots,f_{n+1}$ 
generate $\mathfrak{p}_{n}$. We have proved in Proposition \ref{f1f2} that 
$f_{1},f_{2}\in \mathfrak{p}_{n}$. By our construction, $f_{i}$'s must 
satisfy 
$$a_{(i,i)}zf_{i}+a_{(i,i+1)}yf_{i+1}+a_{(i,i+2)}xf_{i+2}=0,$$ 
with $a_{(i,i+2)}\neq 0$ for $1\leq i\leq n-1$. Therefore 
$f_{i}\in \mathfrak{p}_{n}$, for $1\leq i\leq n+1$.
\medskip

Note that $f_{i}^{\sigma}\in V_{n^{2}+n+i-1} $, for $1\leq i\leq n+1$, where 
\begin{eqnarray*}
V_{r} & = &\{\sigma{\rm -homogeneous \, form \, of} \, \sigma {\rm -order}\,\, r\} \cap \\
{} & {} & \{ \sigma{\rm -leading \, forms \, of \, elements \, in\, \mathfrak{p}_{n}}\}\cup \textbf{0}.
\end{eqnarray*}   
At-first we show that $\{\bar{f}_{1},\ldots,\bar{f}_{n+1}\}$ is linearly independent in 
$\mathfrak{p}_{n}/\mathfrak{m}\mathfrak{p}_{n}$, where $\mathfrak{m} = \langle x,y,z\rangle$ 
is the maximal ideal of the ring $k[[x,y,z]]$.
We consider the expression $\sum_{i=1}^{n+1} h_{i}f_{i}\in\mathfrak{m}\mathfrak{p}_{n}$. 
Let us take a generating set 
$\{p_{1},\ldots,p_{n+1}\}$ of $\mathfrak{p}_{n}$ as described in Moh's 
theorem in section 1. We have $\sum_{i=1}^{n+1} h_{i}f_{i}-\sum_{j=1}^{n+1} d_{j}p_{j}=0$, 
where $d_{j}\in \mathfrak{m}$ for $1\leq j\leq n+1$. Considering the $\sigma$-leading 
form of $\sigma$-order $(n^{2}+n+i-1)$ in the above expression, we get 
$h_{i}\in \mathfrak{m}$. Therefore  $\{\bar{f}_{1},\ldots,\bar{f}_{n+1}\}$ is 
linearly independent in $\mathfrak{p}_{n}/\mathfrak{m}\mathfrak{p}_{n}$. 
Let $\{f_{1},\ldots,f_{n},g\}$ be a minimal generating set of $\mathfrak{p}_{n}$, 
we may assume that this set of generators satisfies the conditions of Moh's theorem 
in section 1. Let $g^{\sigma}$ be the $\sigma$-leading form of $g$. Then by the 
Lemma in section 3 of \cite{moh2}, there is a non trivial relation 
$$b_{(n-1)(n-1)}(0,0,0)zf^{\sigma}_{n-1}+b_{(n-1)n}(0,0,0)yf^{\sigma}_{n}+b_{(n-1)(n+1)}(0,0,0)xg^{\sigma}=0.$$ 
Therefore, by construction of $f_{n+1}^{\sigma}$, we have $g^{\sigma}=cf_{n+1}^{\sigma}$, for some $c\in k$. 
Since $xf_{1}^{\sigma}, g^{\sigma}$ generate $V_{n^{2}+2n}$, we get that 
$xf_{1}^{\sigma}, f_{n+1}^{\sigma}$ also generate $V_{n^{2}+2n}$. Hence 
by Moh's theorem in section 1, $f_{1},\ldots,f_{n+1}$ generates $\mathfrak{p}_{n}$.\qed

\bibliographystyle{amsalpha}

\begin{thebibliography}{A}
\bibitem{coxetal} D. Cox, J. Little, D. O'Shea, \emph{Ideals, Varieties 
and Algorithms}; Springer - Verlag; New York, 1996. 

\bibitem{singular} W. Decker; G.-M. Greuel; G. Pfister; H. Sch\"{o}nemann: Singular 4-1-1 -- A computer algebra system for polynomial computations. http://www.singular.uni-kl.de (2018).

\bibitem{moh1} T. T. Moh, \emph{On the unboundedness of generators of prime ideals in power series rings of three variables}, J. Math. Soc. Japan 26(1974), 722-734.

\bibitem{moh2} T. T. Moh, \emph{On Generators of Ideals}, Proceedings of the American Mathematical Society 77 (3) (1979), 309-312.
\end{thebibliography}

\end{document}